\documentclass{amsart}
\numberwithin{equation}{section}

\newtheorem{thm}{Theorem}[section]

\newtheorem{cor}[thm]{Corollary}

\newtheorem{defin}[thm]{Definition}
\newtheorem{rem}[thm]{Remark}

\begin{document}
\title{Inverse problem of determining an order of the
Riemann-Liouville time-fractional derivative}

\author{Shavkat Alimov}
\author{Ravshan Ashurov}
\address{National University of Uzbekistan named after Mirzo Ulugbek and Institute of Mathematics, Uzbekistan Academy of Science}
\curraddr{Institute of Mathematics, Uzbekistan Academy of Science,
Tashkent, 81 Mirzo Ulugbek str. 100170} \email{ashurovr@gmail.com}

\small

\title[On the backward problems ...] {On the backward problems in time for time-fractional subdiffusion equations}

\begin{abstract}
The backward problem for subdiffusion equation with the fractional
Riemann - Liouville time-derivative of order $\rho\in (0,1)$ and
an arbitrary positive self-adjoint operator $A$ is considered.
This problem is ill-posed in the sense of Hadamard due to the lack
of stability of the solution. Nevertheless, we will show that if
we consider sufficiently smooth current information, then the
solution exists and it is unique. Using this result, we study the
inverse problem of initial value identification for subdiffusion
equation. The results obtained differ significantly from the
corresponding results for the classical diffusion equation (i.e.
$\rho=1$) and even for the subdiffusion equation with the Caputo
derivative. A list of examples of operator $A$ is discussed,
including linear systems of fractional differential equations,
differential models with involution, fractional Sturm-Liouville
operators, and many others.

\end{abstract}

\maketitle

\section{Main results}

The phenomenon of diffusion that does not satisfy the classical
Newton's laws is called anomalous diffusion (subdiffusion or
superdiffusion). Since integer-order diffusion equations cannot
accurately describe anomalous diffusion phenomena in different
fields, and fractional derivatives have the advantages of memory,
they can more accurately describe these anomalous diffusion
phenomena (see, for example, \cite{Koch}-\cite{SU1}). Modeling of
subdiffusion processes is carried out by replacing the first
time-derivative by a fractional one (of the order of $ \alpha \in
(0,1) $) in the classical diffusion equations, and the resulting
equation is called the subdiffusion equation.

When considering subdiffusion equation as model equation in
analyzing anomalous diffusion processes, the data of the model
such as the initial data, the diffusion coefficient, the source
term or even the order of derivative $\alpha$ are not all known,
this leads to study a fractional inverse problem. For this reason,
an additional measurement data is required to deal with this type
of problems. It should be noted, that numerous contributions are
introduced to resolve various fractional inverse problems (see,
for example, \cite{LiLiu}-\cite{Tuan1}).

In this paper, we first consider the backward problem for
subdiffusion equation with Riemann-Liouville fractional derivative
in time and an arbitrary positive selfadjoint operator, having a
discrete spectrum. At the end of the article, we will make the
necessary remarks that allow us to consider self-adjoint operators
with an arbitrary spectrum. It may be worth mentioning that this
problem is ill-posed in the sense of Hadamard due to the lack of
stability of the solution. Nevertheless, we will show that if we
consider sufficiently smooth current information, then the
solution exists and is unique. Using this result, we study the
inverse problem of initial value identification for subdiffusion
equation.

Let us move on to an accurate description of the research objects
and formulate the main results of the work.

Let $H$ be a separable Hilbert space with the scalar product
$(\cdot, \cdot)$ and the norm $||\cdot||$ and  $A: H\rightarrow H$
be an arbitrary unbounded positive selfadjoint operator in $H$.
Suppose that $A$ has a complete in $H$ system of orthonormal
eigenfunctions $\{v_k\}$ and a countable set of nonnegative
eigenvalues $\lambda_k$. It is convenient to assume that the
eigenvalues do not decrease as their number increases, i.e.
$0<\lambda_1\leq\lambda_2 \cdot\cdot\cdot\rightarrow +\infty$.

Using the definitions of a strong integral and a strong
derivative, fractional analogues of integrals and derivatives can
be determined for vector-valued functions (or simply functions)
$h: \mathbb{R}_+\rightarrow H$, while the well-known formulae and
properties are preserved (see, for example, \cite{Liz}). Recall
that the fractional integration of order $ \sigma <0 $ of the
function $ h (t) $ defined on $ [0, \infty) $ has the form
$$
\partial_t^\sigma h(t)=\frac{1}{\Gamma
(-\rho)}\int\limits_0^t\frac{h(\xi)}{(t-\xi)^{\sigma+1}} d\xi,
\quad t>0,
$$
provided the right-hand side exists. Here $\Gamma(\sigma)$ is
Euler's gamma function. Using this definition one can define the
Riemann - Liouville fractional derivative of order $\rho$,
$0<\rho< 1$, as
$$
\partial_t^\rho h(t)= \frac{d}{dt}\partial_t^{\rho-1} h(t).
$$
If in this definition we interchange differentiation and
fractional integration, then we get the definition of the
regularized derivative, that is, the definition of the fractional
derivative in the sense of Caputo:
$$
D_t^\rho h(t)= \partial_t^{\rho-1}\frac{d}{dt} h(t).
$$

Note that if $\rho=1$, then fractional derivatives coincides with
the ordinary classical derivative of the first order: $\partial_t
h(t) = D_t h(t)= \frac{d}{dt} h(t)$.

Let $\rho\in(0,1) $ be a fixed number and let $C((a,b); H)$ stand
for a set of continuous functions $u(t)$  of $t\in (a,b)$ with
values in $H$. The space $C^\infty((a,b); H)$ is defined
similarly.

Consider the Cauchy type problem with inverse time:
\begin{equation}\label{prob1}
\left\{
\begin{aligned}
&\partial_t^\rho u(t) + Au(t) = f(t),\quad 0<t< T;\\
&u(T) = \Phi,
\end{aligned}
\right.
\end{equation}
where $\Phi\in H$ and $f(t)\in C((0,T); H)$ are  given vectors.
This problem is called \emph{the backward problem} (see, for
example, \cite{Yama10}-\cite{Florida}).

\begin{defin}\label{def} A function $u(t)\in C((0,T]; H)$ with the properties
$\partial_t^\rho u(t), Au(t)\in C((0,T); H)$ and satisfying
conditions (\ref{prob1})  is called \textbf{the
 solution} of the backward problem (\ref{prob1}).
\end{defin}

The standard formulation of the Cauchy problem for equation
(\ref{prob1}) has the form:
\begin{equation}\label{prob2}
\left\{
\begin{aligned}
&\partial_t^\rho u(t) + Au(t) = f(t),\quad 0<t< T;\\
&\lim\limits_{t\rightarrow 0}\partial_t^{\rho-1}u(t)=\varphi,
\end{aligned}
\right.
\end{equation}
where $\varphi\in H$ is a given vector. This problem will be
called \emph{the forward problem}. The solution to this problem is
defined similarly to the solution to the backward problem. In
order to investigate the backward problem, one usually uses the
properties of the solution to the forward problem.

Let $ \tau $ be an arbitrary real number. We introduce the power
of operator $ A $, acting in $ H $ according to the rule
$$
A^\tau g= \sum\limits_{k=1}^\infty \lambda_k^\tau g_k v_k,
$$
where $g_k$ is the Fourier coefficients of a function $g\in H$:
$g_k=(g, v_k)$. Obviously, the domain of this operator has the
form
$$
D(A^\tau)=\{g\in H:  \sum\limits_{k=1}^\infty \lambda_k^{2\tau}
|g_k|^2 < \infty\}.
$$
For elements of $D(A^\tau)$ we introduce the norm
\[
||g||^2_\tau=\sum\limits_{k=1}^\infty \lambda_k^{2\tau} |g_k|^2 =
||A^\tau g||^2.
\]

We first prove the existence and uniqueness of a solution of
problem (\ref{prob1}).

\begin{thm}\label{mt} Let $f(t)\equiv 0$. Then for any $\Phi\in D(A^2)$ problem (\ref{prob1}) has a unique
solution. Moreover there exist constants $C_1, C_2>0$, such that
\begin{equation}\label{bph}
C_1\lim\limits_{t\rightarrow 0}||\partial_t^{\rho-1}u(t)||\leq
||u(T)||_2\leq C_2\lim\limits_{t\rightarrow
0}||\partial_t^{\rho-1}u(t)||.
\end{equation}
\end{thm}

\begin{rem}If, for example, $A$ is an elliptic operator of the second order,
then in order for a solution to the backward problem (\ref{prob1})
to exist, $\Psi$ must have four derivatives.

\end{rem}

The backward problems for the diffusion process are of great
importance in engineering fields and are aimed at determining the
previous state of a physical field (for example, at $t=0$) based
on its current information (see, for example, \cite{Yama11}).
However, regardless of the fact that the Riemann-Liouville or
Caputo derivative is taken into the equation, this problem is
ill-possed in the sense of Hadamard. In other words, a small
change of $u(T)$ in the norm of space $H$ leads to large changes
in the initial data. As can be seen from the above theorem, the
situation changes if we take the norm of space $D(A^2)$ instead of
the norm in $H$. It should also be noted, that for the backward
problem of the classical diffusion equation (that is $\rho =1$)
estimates of the type (\ref{bph}) on the scales of spaces $D(A^a)$
are generally impossible (see, for example, Chapter 8.2 of
\cite{Kab1}).

In case of the Caputo derivative $D_t^\rho$ the problem
(\ref{prob1}) for various elliptic differential operators $A$ has
been considered by a number of authors. Let us mention only some
of these works. For the case of the second order symmetric
elliptic operator $A$, Sakamoto and Yamamoto \cite{Yama11}
establish the unique existence of weak solutions and the
asymptotic behavior as time $t$ goes to $\infty$. They also prove
the stability in the backward problem in time and the uniqueness
in determining an initial value. Nonsymmetric case was considered
in Florida, Li, Yamamoto \cite{Florida}. Since the problem is
ill-possed, many authors have considered various regularization
options for finding the initial condition (see, for
one-dimensional elliptical part, Liu and Yamamoto \cite{Yama10},
for the nonlinear case,  Tuan, Huynh, Ngoc, and Zhou
\cite{Tuan1}). In particular, as for numerical approaches, see
Tuan, Long and Tatar \cite{Tuan2}, Wang and Liu \cite{Wang} and
the references  therein.

For backward problem (\ref{prob1})  with non-homogeneous term we
have the following result (for the Caputo derivative $D_t^\rho$
see the above mentioned work \cite{Florida}):

\begin{thm}\label{tbpn} Let $t^{1-\rho}f(t)\in C([0,T];
D(A^{1+\varepsilon}))$ with some $\varepsilon>0$. Then for any
$\Phi\in D(A^2)$ problem (\ref{prob1}) has a unique solution.
Moreover there exists a constant $C>0$, such that
\begin{equation}\label{bpn}
\lim\limits_{t\rightarrow 0}||\partial_t^{\rho-1}u(t)||\leq C\big(
||u(T)||_2+\max\limits_{t\in[0,T]}||t^{1-\rho}f(t)||_{1+\varepsilon}\big).
\end{equation}

\end{thm}

The remainder of this paper is composed of three sections. Section
2 is devoted to the study of forward problem (\ref{prob2}). In
Section 3, we show that problem (\ref{prob1}) is ill-possed in the
sense of Hadamard and prove Theorems \ref{mt} and \ref{tbpn}. In
the last section, examples of the operator $A$ are presented. In
addition, the necessary remarks are given, with the help of which
all the statements formulated can be translated to the case when
the operator $A$ has an arbitrary spectrum.

Finally we note, that to investigate the forward and backward
problems we borrow some original ideas from papers \cite{Yama11},
\cite{Florida}, where authors studied the similar problems for
equation (\ref{prob1}) with the Caputo derivative.

\section{Forward problem}

\begin{thm}\label{fpn} Let $\varphi\in H$ and $t^{1-\rho}f(t)\in
C([0,T]; D(A^\varepsilon))$ for some $\varepsilon\in (0, 1)$. Then
there exists a unique solution to the forward problem, such that
\begin{equation}\label{fpn}
\left\{
\begin{aligned}
&||u(t)||_1+||\partial^\rho_t u(t)||\leq C_\varepsilon\big(
t^{-1-\rho}||A^{-1}\varphi||+\max\limits_{t\in[0,T]}||t^{1-\rho}f(t)||_{\varepsilon}+||f(t)||\big), \quad t>0,\\
&\max\limits_{t\in[0,T]}||t^{1-\rho}u(t)||\leq C
(||\varphi||+\max\limits_{t\in[0,T]}||t^{1-\rho}f(t)||),
\end{aligned}
\right.
\end{equation}
where $ C $ is an absolute constant and $ C_\varepsilon $ is a
constant depending on $ \varepsilon $.

Moreover, if $t^{1-\rho}f(t)\in C([0,T]; D(A^{1+\varepsilon}))$,
then there exists a constant $C$, depending on $\varepsilon$ such
that
\begin{equation}\label{fpn0}
||u(t)||_2\leq C(
t^{-1-\rho}||\varphi||+\max\limits_{t\in[0,T]}||t^{1-\rho}f(t)||_{1+\varepsilon}),
\quad t>0.
\end{equation}
\end{thm}

Here is an obvious consequence of estimate (\ref{fpn0}):
\begin{cor}\label{C} Let $\varphi\in H$ and $t^{1-\rho}f(t)\in C([0,T];
D(A^{1+\varepsilon}))$. Then there exists a constant $C$,
depending on $T$ and $\varepsilon$ such that
\begin{equation}\label{corE}
||u(T)||_2\leq C(
||\varphi||+\max\limits_{t\in[0,T]}||t^{1-\rho}f(t)||_{1+\varepsilon}).
\end{equation}

\end{cor}

\begin{thm}\label{fpHD} Let $\varphi\in H$ and $f\equiv 0$. The
the unique solution to the forward problem is infinitely
differentiable with respect to the variable $t$, i.e.
\[
u(t)\in C^\infty ((0,\infty); H),
\]
and there exists a constant $C$ such that the following estimates
are valid
\begin{equation}\label{fphE2}
\left\{
\begin{aligned}
&||u(t)||\leq \frac{C t^{\rho-1}}{1+(\lambda_1 t^\rho)^2}\cdot||\varphi||, \quad t> 0,\\
&||\partial^m_t u(t)||\leq C t^{\rho-1-m}\cdot ||\varphi||, \quad
t>0, \quad m\in \mathbb{N}.
\end{aligned}
\right.
\end{equation}

\end{thm}

We note at once that the solution to equation (\ref{prob1}),
generally speaking, is not continuous at the point $t = 0$ (see
(\ref{fpn})). Of course, we can consider a continuous function at
the point $t = 0$ as the right-hand side of the equation, but we
assume that $f(t)$ would have a singularity at this point in order
to cover a more general case (see (\ref{fpn}) and (\ref{fpn0})).

Initial-boundary value problems for various subdiffusion equations
have been investigated by many specialists. Let us mention only
some of these works. In the book of A.A. Kilbas et al. \cite{Kil}
(Chapter 6) there is a survey of works published before 2006. The
case of one spatial variable $x\in \mathbb{R}$ and subdiffusion
equation with "the elliptical part" $u_{xx}$ were considered  for
example in the monograph of A. V. Pskhu \cite{PSK} (Chapter 4, see
references thesein). The paper Gorenflo, Luchko and Yamamoto
\cite{GorLuchYam} is devoted to the study of subdiffusion
equations in Sobelev spaces. In the paper by Kubica and Yamamoto
\cite{KubYam}, initial-boundary value problems for equations with
time-dependent coefficients are considered. In the
multidimensional case ($x\in \mathbb{R}^N$), instead of the
differential expression $u_{xx}$, authors considered either the
second order elliptic operator (\cite{Agr} - \cite{PS1}) or
elliptic pseudodifferential operators with constant coefficients
in the whole space $\mathbb{R}^N$ (Umarov \cite{SU}). In the paper
of Yu. Luchko \cite{Luch} the author constructed solutions by the
eigenfunction expansion in the case of $f = 0$ and discussed the
unique existence of the generalized solution to problem
(\ref{prob2}) with the Caputo derivative. The authors of the
recent paper \cite{AO} considered initial-boundary value problems
for subdiffusion equations with arbitrary elliptic differential
operators in bounded domains.

The formulated results for equation (\ref{prob1}) with the Caputo
derivative were previously proved in \cite{Yama11} (in the case
when $A$ is a symmetric elliptic operator of the second order) and
\cite{Florida} (in the case when $A$ is not symmetric).

Let us compare our results (equation (\ref{prob1}) with the
Riemann-Liouville derivative) with the results of \cite{Yama11}
and \cite{Florida} (equation (\ref{prob1}) with the Caputo
derivative) and standard results for the case of $\rho=1$.

In our case we have no smoothing properties like the case of the
Caputo and the classical diffusion equation (i.e. $\rho=1$). In
Theorem \ref{fpn} (estimate (\ref{fpn0}) there is the smoothing
property in space with order 2 which means that $u(t)\in D(A^2)$
for any $t > 0$ and any $\varphi\in H$. For example, if $A$ is an
elliptic operator of order two, defined in $N$-dimensional domain
$\Omega$, then the condition $\varphi\in L_2(\Omega)$ guarantee
that the solution to problem (\ref{prob2}) is in the classical
Sobolev space $W_2^4(\Omega)$. Nevertheless, as it is proved in
Theorem \ref{fpHD} the regularity in time immediately becomes
stronger in $t$, and is of infinity order (i.e., $u$ is of
$C^\infty$ for $t > 0$). In Theorem \ref{mt}, it is  showed that
the smoothing in $D(A^2)$ is the best possible and the solution
cannot be smoother than $D(A^2)$ at $t > 0$ if $\varphi\in H$.

For the case of the Caputo derivative, in papers \cite{Yama11} and
\cite{Florida} it is proved that the best possible smoothing
property is of order 1, i.e. $u(t)\in D(A)$ for any $t > 0$ and
any $\varphi\in H$, while in the classical case ($\rho=1$) the
solution is infinitely differentiable both with respect to spatial
variables and $t$ with any $t > 0$ and $\varphi\in H$.

The first estimate (\ref{fphE2}) shows the decay of solution with
order $t^{-\rho-1}$ as $t\rightarrow\infty$, which is slower than
the exponential decay in the case of $\rho=1$, but faster than the
Caputo case with the decay of $t^{-\rho}$ (see \cite{Yama11}).

\textbf{Proof of Theorem \ref{fpn}.} Let us introduce the
following formal series
\begin{equation}\label{sp1}
u(t)= \sum\limits_{k=1}^\infty \big[t^{\rho-1} E_{\rho,
\rho}(-\lambda_k t^{\rho}) \varphi_k + \int\limits_0^t
\eta^{\rho-1} E_{\rho, \rho} (-\lambda_k \eta^\rho)f_k(t-\eta)
d\eta \big]v_k,
\end{equation}
where $\varphi_k$ and $f_k(t)$ are the Fourier coefficients of
$\varphi$ and $f(t)$ correspondingly, $E_{\rho, \mu}(t)$ - the
Mittag-Leffler function:
$$
E_{\rho, \mu}(t)= \sum\limits_{n=0}^\infty \frac{t^n}{\Gamma(\rho
n+\mu)}.
$$

By virtue of the formula (\cite{PSK}, p. 104)
\begin{equation} \label{property}
\lim\limits_{t\rightarrow +0}\partial_t^{\alpha-1} u(t) = \Gamma
(\alpha)\lim\limits_{t\rightarrow +0} t^{1-\alpha} u(t).
\end{equation}
one can easily verify that the function (\ref{sp1}) formally
satisfies the conditions of problem (\ref{prob2}) (see, for
example, \cite{Gor}, p. 173). In order to prove that function
(\ref{sp1}) is actually a solution to the problem, it remains to
substantiate this formal statement, i.e. show that the operators
$A$ and $\partial_t^\rho$ can be applied term by term to the
series (\ref{sp1}).

To do this we need the asymptotic estimate of the Mittag-Leffler
function with a sufficiently large  negative argument. The well
known estimate has the form (see, for example, \cite{Dzh66}, p.
136)
\[
|E_{\rho, \mu}(-t)|\leq \frac{C}{1+ t}, \quad t>0,
\]
where $\mu$ is an arbitrary complex number. But for $E_{\rho,
\rho}(-t)$ one can get a better estimate. Indeed, using the
asymptotic estimate (see, for example, \cite{Dzh66}, p. 134)
\begin{equation}\label{mA}
E_{\rho, \rho}(-t)=-\frac{t^{-2}}{\Gamma(-\rho)} +O(t^{-3}).
\end{equation}
and the fact that $E_{\rho, \rho}(t)$ is real analytic, we can
obtain the following inequality
\begin{equation}\label{m1}
|E_{\rho, \rho}(-t)|\leq \frac{C}{1+ t^2}, \quad t>0.
\end{equation}

We will also use a coarser estimate with positive eigenvalues
$\lambda_k$ and $0<\varepsilon<1$:
\begin{equation}\label{m2}
|t^{\rho-1} E_{\rho,\rho}(-\lambda_k t^\rho)|\leq
\frac{Ct^{\rho-1}}{1+(\lambda_kt^\rho)^2}\leq C
\lambda^{\varepsilon-1}_k t^{\varepsilon\rho-1}, \quad t>0,
\end{equation}
which is easy to verify. Indeed, let $t^\rho\lambda_k<1$, then $t<
\lambda_k^{-1/\rho}$ and
$$
t^{\rho -1} = t^{\rho-\varepsilon\rho} t^{\varepsilon\rho-1} <
\lambda_k^{\varepsilon-1}t^{\varepsilon\rho-1}.
$$
If $t^\rho\lambda_k\geq 1$, then $\lambda_k^{-1}\leq t^\rho$ and
$$
\lambda_k^{-2} t^{-\rho-1}=\lambda_k^{-1+\varepsilon}
\lambda_k^{-1-\varepsilon} t^{-\rho-1}\leq
\lambda_k^{\varepsilon-1}t^{\varepsilon\rho-1}.
$$

Let $S_j(t)$ be the partial sum of series (\ref{sp1}). Then
$$
A S_j(t)=\sum\limits_{k=1}^j \big[t^{\rho-1} E_{\rho,
\rho}(-\lambda_k t^{\rho}) \varphi_k + \int\limits_0^t
\eta^{\rho-1} E_{\rho, \rho} (-\lambda_k \eta^\rho)f_k(t-\eta)
d\eta \big]\lambda_k v_k.
$$
Due to the Parseval equality we may write
$$
||AS_j(t)||^2=\sum\limits_{k=1}^j \lambda^2_k \big|t^{\rho-1}
E_{\rho, \rho}(-\lambda_k t^\rho)\, \varphi_k+\int\limits_0^t
\eta^{\rho-1} E_{\rho, \rho} (-\lambda_k \eta^\rho)f_k(t-\eta)
d\eta\big|^2.
$$
Using estimate (\ref{m1}) and the inequality $\lambda_j
t^{\rho-1}(1+(\lambda_k t^{\rho})^2)^{-1} < \lambda_k^{-1}
t^{-\rho-1}$ we obtain
$$
 \sum\limits_{k=1}^j \lambda^2_k \big|t^{\rho-1}
E_{\rho, \rho}(-\lambda_k t^\rho)\, \varphi_k|^2\leq C
t^{-2(1+\rho)} \sum\limits_{k=1}^j \lambda^{-2}_k  |\varphi_k|^2=
C t^{-2(1+\rho)}||A^{-1}\varphi||^2.
$$
On the other hand, by inequality (\ref{m2}) for $0<\varepsilon<1$
one has
\[
\sum\limits_{k=1}^j \lambda^2_k \bigg|\int\limits_0^t
\eta^{\rho-1} E_{\rho, \rho} (-\lambda_k \eta^\rho)f_k(t-\eta)
d\eta\bigg|^2\leq C\sum\limits_{k=1}^j \bigg[ \int\limits_0^t
\eta^{\varepsilon\rho-1}\lambda_k^{\varepsilon}|f_k(t-\eta)|
d\eta\bigg]^2\leq
\]
(by virtue of the generalized Minkowski inequality)
\[
\leq
C\bigg[\int\limits_0^t\eta^{\varepsilon\rho-1}(t-\eta)^{\rho-1}\bigg(\sum\limits_{k=1}^j
|\lambda_k^{\varepsilon}(t-\eta)^{(1-\rho)}f_k(t-\eta)|^2\bigg)^{\frac{1}{2}}
d\eta\bigg]^2\leq C_\varepsilon
\max\limits_{t\in[0,T]}||t^{1-\rho}f(t)||^2_{\varepsilon}.
\]
Therefore
\begin{equation}\label{Aj}
||AS_j(t)||^2\leq C
t^{-2(1+\rho)}||A^{-1}\varphi||^2+C_\varepsilon
\max\limits_{t\in[0,T]}||t^{1-\rho}f(t)||^2_\varepsilon, \quad
t>0.
\end{equation}
Hence, we obtain $Au(t)\in C((0,T]; H)$.

Further, from equation (\ref{prob1}) one has $\partial_t^\rho
S_j(t)= - AS_j(t)+\sum\limits_{k=1}^jf_k(t) v_k$. Therefore, from
above reasoning, we have $\partial_t^\rho u(t)\in C((0,T]; H)$ and
\begin{equation}\label{Sj}
||\partial_t^\rho S_j(t)||^2\leq C
t^{-2(1+\rho)}||A^{-1}\varphi||^2+C_\varepsilon
\max\limits_{t\in[0,T]}||t^{1-\rho}f(t)||^2_{\varepsilon}+||f(t)||^2,
\quad t>0.
\end{equation}
Thus, we have completed the rationale that (\ref{sp1}) is a
solution to the forward problem. Inequalities (\ref{Aj}) and
(\ref{Sj}) imply the first estimate (\ref{fpn}). The second
estimate (\ref{fpn}) is proved using similar reasoning.

Now let $t^{1-\rho}f(t)\in C([0,T]; D(A^{1+\varepsilon}))$. Then
using estimate (\ref{m1}) and the inequality $\lambda^2_j
t^{\rho-1}(1+(\lambda_k t^{\rho})^2)^{-1} < t^{-\rho-1}$ we obtain
$$
\sum\limits_{k=1}^j \lambda^4_k \big|t^{\rho-1} E_{\rho,
\rho}(-\lambda_k t^\rho)\, \varphi_k\big|^2\leq C
t^{-2(\rho+1)}||\varphi||, \quad t>0.
$$
On the other hand, by inequality (\ref{m2}) we have
\[
\sum\limits_{k=1}^j \lambda^4_k \bigg|\int\limits_0^t
\eta^{\rho-1} E_{\rho, \rho} (-\lambda_k \eta^\rho)f_k(t-\eta)
d\eta\bigg|^2\leq C\sum\limits_{k=1}^j \bigg[ \int\limits_0^t
\eta^{\varepsilon\rho-1}\lambda_k^{1+\varepsilon}|f_k(t-\eta)|
d\eta\bigg]^2\leq
\]
\[
\leq C_\varepsilon
\max\limits_{t\in[0,T]}||t^{1-\rho}f(t)||^2_{1+\varepsilon}.
\]
Therefore
\[
||A^2 S_j(t)||\leq C t^{-2(\rho+1)}||\varphi|| + C_\varepsilon
\max\limits_{t\in[0,T]}||t^{1-\rho}f(t)||^2_{1+\varepsilon}, \quad
t>0.
\]
This implies estimate (\ref{fpn0}).

We now turn to the proof of the the uniqueness of the forward
problem's solution.

Suppose that problem (\ref{prob1}) has two solutions $u_1(t)$ and
$u_2(t)$. Our aim is to prove that $u(t)=u_1(t)-u_2(t)\equiv 0$.
Since the problem is linear, then we have the following homogenous
problem for $u(t)$:
\begin{equation}\label{ur1}
\partial_t^\rho u(t) + Au(t) = 0, \quad t>0;
\end{equation}
\begin{equation}\label{nu1}
\lim\limits_{t\rightarrow 0}\partial_t^{\rho-1} u(t) = 0.
\end{equation}

Set
\[
w_k(t)\ =\ (u(t), v_k).
\]
It follows from (\ref{ur1}) that for any $k\in \mathbb{N}$
\[
\partial_t^\rho w_k(t) = (\partial_t^\rho u(t), v_k) = -\, (Au(t), v_k)= -\,
(u(t), Av_k) = -\, \lambda_k w_k(t).
\]
Therefore, we have the following Cauchy problem for $w_k(t)$ (see
(\ref{nu1})):
$$
\partial_t^\rho w_k(t) +\lambda_k w_k(t)=0,\quad t>0; \quad \lim\limits_{t\rightarrow 0}\partial_t^{\rho-1} w_k(t)=0.
$$
This problem has the unique solution (see, for example,
\cite{Gor}, p. 173 and \cite{ACT}). Therefore, $w_k(t) = 0$ for
$t>0$ and for all $k\geq 1$. Then by the Parseval equation we
obtain $u(t) = 0$ for all $t>0$. Hence uniqueness of the solution
is proved.

Thus the proof of Theorem \ref{fpn} is complete.

\

\textbf{Proof of Theorem \ref{fpHD}.} It is sufficient to prove
the estimates (\ref{fphE2}) since the infinitely differentiability
of the solution follows from the second estimate (\ref{fphE2}).

The estimate of the Mittag-Leffler function (\ref{m1}) implies for
$t>0$ the first inequality (\ref{fphE2}):
\[
||u(t)||^2= \sum\limits_{k=1}^\infty \big[t^{\rho-1} E_{\rho,
\rho}(-\lambda_k t^{\rho}) \varphi_k]^2\leq
\sum\limits_{k=1}^\infty |\varphi_k|^2 \bigg(\frac{C
t^{\rho-1}}{1+ (\lambda_k t^\rho)^2}\bigg)^2\leq \bigg(\frac{C
t^{\rho-1}}{1+ (\lambda_1 t^\rho)^2} ||\varphi||\bigg)^2.
\]
To prove the second estimate we remind the following
differentiation formula for the Mittag-Leffler function (see, for
example, \cite{Gor}, formula (4.3.1))
\[
\bigg(\frac{d}{dz}\bigg)^m\bigg[z^{\mu-1} E_{\alpha, \mu}
(z^\alpha)\bigg]=z^{\mu-1-m}E_{\alpha, \mu-m} (z^\alpha),\quad
z\neq 0, \quad m\in \mathbb{N},
\]
which is an immediate consequence of the definition of the
Mittag-Leffler function $E_{\alpha, \mu}$. Therefore one has
\[
\bigg(\frac{d}{dt}\bigg)^m\bigg[t^{\rho-1} E_{\rho, \rho}
(-\lambda_k t^\rho)\bigg]=t^{\rho-1-m}E_{\rho, \rho-m} (-\lambda_k
t^\rho),\quad t>0, \quad m\in \mathbb{N},
\]
By this formula, we have (note $\partial_t^m = d^m/dt^m$)
\[
\partial_t^m u(t)= \sum\limits_{k=1}^\infty t^{\rho-1-m} E_{\rho,
\rho-m}(-\lambda_k t^{\rho}) \varphi_k v_k, \quad t>0,
\]
for $m\in \mathbb{N}$ and estimate (\ref{m1}) implies
\[
||\partial^m_t u(t)||\leq C t^{\rho-1-m}\cdot ||\varphi||, \quad
t>0, \quad m\in \mathbb{N}.
\]
The proof of Theorem \ref{fpHD} is completed.

\section{Backward problem}

In the case $ \rho = 1 $ problem (\ref{prob1}) is called (see, for
example, \cite{Kab1}, p. 214) the inverse heat conduction problem
with inverse time (retrospective inverse problem). This problem is
ill-posed. Usually, when investigating such inverse problems, they
are reduced by changing the variables $\tau=T-t$ to an equivalent
forward problem, where the sign changes at the time derivative
(\cite{Kab1}, p. 214). However, this approach cannot be
implemented for problem (\ref{prob1}) since the simple property of
derivatives $\frac{d}{dt}=-\frac{d}{d\tau}$ does not hold for
fractional derivatives: $\partial_t^\rho \neq-
\partial_\tau^\rho$.

Problem (\ref{prob1}) is also ill-posed in the sense of Hadamard
because of the same reason as the classical one ($ \rho = 1 $): a
small variation of $u(T)$ in the norm of space $H$ may cause
arbitrarily large variations in the initial data. Indeed, let
$f(t)\equiv 0$ and $u(T)=\lambda^{-2+\varepsilon}_k v_k $,
$\varepsilon>0$, in problem (\ref{prob1}). Then the unique
solution of the problem is
\[
u(t)=\lambda^{-2+\varepsilon}_k\cdot\frac{t^{\rho-1}E_{\rho,
\rho}(-\lambda_k t^\rho)}{T^{\rho-1}E_{\rho, \rho}(-\lambda_k
T^\rho)} \cdot v_k
\]
and
\[
\lim\limits_{t\rightarrow +0}\partial_t^{\rho-1}
u(t)=\lambda^{-2+\varepsilon}_k\cdot\frac{1}{T^{\rho-1}E_{\rho,
\rho}(-\lambda_k T^\rho)} \cdot v_k.
\]
Therefore, on the one hand, $||u(T)||=\lambda^{-2+\varepsilon}_k$
and it tends to zero as $k\rightarrow \infty$ (even
$||u(T)||_{a}\rightarrow 0$ for any $a<2-\varepsilon$), and on the
other,  according to the asymptotic estimate (\ref{mA}),
\[
||\lim\limits_{t\rightarrow +0}\partial_t^{\rho-1}
u(t)||=\lambda^{-2+\varepsilon}_k\cdot\frac{1}{T^{\rho-1}E_{\rho,
\rho}(-\lambda_k T^\rho)} \rightarrow \infty \quad \text{when}
\quad k\rightarrow\infty.
\]

However, if we consider the norm of $u(T)$ in space $D(A^2)$, then
the situation will change completely; note the norm
$||u(T)||_2=\lambda_k^\varepsilon$ in this example is unbounded as
$k\rightarrow\infty$.

We now turn to the proof of Theorems \ref{mt} and \ref{tbpn}.

\textbf{Proof of Theorem \ref{mt}.} Since function $E_{\rho,
\rho}(z)$ has no negative zero (see, for example, \cite{Gor}, p.
74) and $E_{\alpha, \alpha}(0)= \Gamma^{-1}(\rho)>0$, then
\begin{equation}\label{E0}
E_{\rho, \rho}(-t)>0, \quad t\geq 0.
\end{equation}
Let $\Phi\in D(A^2)$ and $\Phi_k$ be its Fourier coefficients.
Then
\[
||\Phi||_2=\sum\limits_{k=1}^\infty \lambda_k^4|\Phi_k|^2<\infty.
\]
By (\ref{E0}) we can set
\[
\varphi_k=\frac{\Phi_k}{T^{\rho-1} E_{\rho, \rho}(-\lambda_k
T^\rho)}.
\]
Then by virtue of the asymptotic estimate (\ref{mA}) one has
\[
\sum\limits_{k=1}^\infty \varphi^2_k=\sum\limits_{k=1}^\infty
\frac{\Phi^2_k}{(T^{\rho-1} E_{\rho, \rho}(-\lambda_k
T^\rho))^2}=\sum\limits_{k=1}^\infty T^{2\rho+2} \lambda_k^4
\Gamma^2(-\rho)\Phi^2_k \bigg(\frac{1}{1+O(\lambda^{-3}_k
T^{-3\rho})}\bigg)^2\leq
\]
\[
\leq CT^{2\rho+2}\sum\limits_{k=1}^\infty
\lambda_k^4|\Phi_k|^2<\infty.
\]
Therefore
\[
\varphi=\sum\limits_{k=1}^\infty \varphi_k v_k \in H,
\]
and the following function (see (\ref{sp1}))
\[
u(t)= \sum\limits_{k=1}^\infty t^{\rho-1} E_{\rho,
\rho}(-\lambda_k t^{\rho}) \varphi_k v_k
\]
is the unique solution to forward problem (\ref{prob2}) with
$f(t)\equiv 0$ and the initial function $\varphi$. Moreover
$u(T)=\Phi$ and
\[
\lim\limits_{t\rightarrow 0} ||t^{1-\rho} u(t)||=||\varphi||\leq
C||\Phi||_2= C ||u(T)||_2.
\]
The second inequality in (\ref{bph}) is already proved in Theorem
\ref{fpn} (estimate (\ref{fpn0})).

Theorem \ref{mt} is proved.

\textbf{Proof of Theorem \ref{tbpn}.} Consider the following two
auxiliary problems (see \cite{Florida}):
\begin{equation}\label{prob3}
\left\{
\begin{aligned}
&\partial_t^\rho v(t) + Av(t) = f(t),\quad 0<t< T;\\
&\lim\limits_{t\rightarrow 0} \partial_t^\rho v(t)=0,
\end{aligned}
\right.
\end{equation}
and
\begin{equation}\label{prob4}
\left\{
\begin{aligned}
&\partial_t^\rho w(t) + Aw(t) = 0,\quad 0<t< T;\\
&w(T) = \Phi-v(T),
\end{aligned}
\right.
\end{equation}
If $t^{1-\rho}f(t)\in C([0,T]; D(A^{1+\varepsilon}))$, then there
exists the unique solution to problem (\ref{prob3}) and (see
Corollary \ref{C})
\begin{equation}\label{v}
||v(T)||_2\leq
C\max\limits_{t\in[0,T]}||t^{1-\rho}f(t)||_{1+\varepsilon}.
\end{equation}
If $\Phi\in D(A^{2})$, then there exists the unique solution to
problem (\ref{prob4}) and (see Theorem \ref{mt})
\begin{equation}\label{w}
\lim\limits_{t\rightarrow 0}||\partial_t^{\rho-1}w(t)||\leq C
||w(T)||_2.
\end{equation}

Setting $u=v+w$, we see that $u(T)=\Phi-v(T)+v(T)=\Phi$. Then we
can verify that $u$ is the unique solution to problem
(\ref{prob1}) and the estimates (\ref{v}) and (\ref{w}) imply
\[
\lim\limits_{t\rightarrow
0}||\partial_t^{\rho-1}u(t)||=\lim\limits_{t\rightarrow
0}||\partial_t^{\rho-1}w(t)||\leq ||w(T)||_2\leq
C(||\Phi||_2+||v(T)||_2)\leq
\]
\[
\leq C(||\Phi||_2
+\max\limits_{t\in[0,T]}||t^{1-\rho}f(t)||_{1+\varepsilon}).
\]
Theorem \ref{tbpn} is proved.

\section{Examples of operator $A$ and further generalization}

The setting of an abstract operator $A$ as in this paper allows
one to include many models. For example, as a example one may
consider any of physical examples, discussed in Section 6 of the
paper of M. Ruzhansky et al. \cite{Ruz}, including Sturm-Liouville
problems, differential models with involution, fractional
Sturm-Liouville operators, harmonic and anharmonic oscillators,
Landau Hamiltonians, fractional Laplacians, and harmonic and
anharmonic operators on the Heisenberg group. It should be noted,
that the authors of \cite{Ruz} considered a class of inverse
problems for restoring the right-hand side of a subdiffusion
equation with the Caputo derivatives of order $0<\rho\leq 1$ for a
large class of positive operators with discrete spectrum.

Usually, when studying the subdiffusion equation, an elliptic
equation of order two on a $N$-dimensional bounded domain $\Omega$
with classical boundary conditions, such as Dirichlet, is
considered as the elliptic part. The system of eigenfunctions of
such operators constitutes a complete set in $L_2(\Omega)$, and
the spectrum is discrete and rather regular, that is,
$N(\lambda)=\sum_{\lambda_k\leq \lambda} 1$ - the number of
eigenvalues not exceeding $\lambda$ has the estimate
$N(\lambda)=O(\lambda^{N/2})$. However, for example, for the
Laplace operator in a bounded domain $\Omega$, boundary conditions
can be specified such that the system of eigenfunctions remains
complete in $L_2(\Omega)$, but the set $\{\lambda_k\}$ will be
dense in $(1, + \infty)$ (see \cite{IF}). It should be noted that
the theorems formulated above are also valid for such operators.

On the other hand, not all operators important for applications
have a discrete spectrum. For example, if $A$ is the Laplace
operator $-\triangle$ in $H=L_2(\mathbb{R}^N)$, then the spectrum
of this operator is continuous. A natural question arises: is it
possible to apply the above reasoning to the case of the operator
$A$ with continuous spectrum? The answer to this question is yes
and similar theorems are true as above. Moreover, there is no need
to make significant changes to the corresponding proofs, except
for the proof of the uniqueness of the solution to the forward
problem. Below we give a proof of uniqueness in the general case.

Let $H$ be a Hilbert space with the scalar product $(\cdot,
\cdot)$ and the norm $||\cdot||$ and  $A: H\rightarrow H$ be an
arbitrary semibounded (with the bound $\mu>0$)  selfadjoint
operator in $H$. By von Neumann's spectral theorem, the operator
$A$ has a partition $\{P_\lambda\}$ of unity, and can be
represented in the form of
\[
A\ =\ \int\limits_\mu^\infty \lambda\, dP_\lambda, \quad \mu>0.
\]
The projections $P_\lambda$ increase monotonically, are continuous
on the left, and tend strongly to the unit operator, that is,
\[
\lim\limits_{\lambda \rightarrow \infty} ||P_\lambda g-g||=0,
\quad g\in H.
\]

For any real number $\tau$ the power of operator $A$ is defined as
\[
A^\tau g = \int\limits_\mu^\infty \lambda^\tau\, dP_\lambda
g,\quad g\in D(A^\tau)=\{g\in H  :  \int\limits_\mu^\infty
\lambda^{2\tau}\, (dP_\lambda g, g)<\infty\}.
\]

Consider the  forward problem (\ref{prob2}) with this operator
$A$. Let us define the space $D(A^\varepsilon)$ as above and let
$\varphi\in H$ and $t^{1-\rho}f(t)\in C([0,T]; D(A^\varepsilon))$
for some $\varepsilon\in (0, 1)$. Then it is not hard to verify
that the following function
\[
u(t)=\int\limits_\mu^\infty t^{\rho-1} E_{\rho, \rho} (-\lambda \,
t^\rho)\, dP_\lambda\varphi +
\int\limits_0^t\int\limits_\mu^\infty \eta^{\rho-1} E_{\rho, \rho}
(-\lambda \,\eta^\rho)\, dP_\lambda f(t-\eta)d \eta
\]
is the solution to the forward problem (\ref{prob2}).

Next, we will show that problem (\ref{prob2}) has a unique
solution.

Let $U$ be the spectral representation of $H$ on the direct sum
$\bigoplus\limits_{\beta\in B} L_2(\mathbb{R}_\mu,\mu_\beta)$
($\mathbb{R}_\mu =(\mu, +\infty)$) with respect to self-adjoint
operator $A$ (see \cite{Dun}, Chapter XII, Sec. 3, Theorem 5).
Operator $U$ is a linear map of $H$ onto all space
$\bigoplus\limits_{\beta\in B} L_2(\mathbb{R}_\mu,\mu_\beta)$ and
preserves the scalar product, that is, it is a unitary operator.

Note that
\[
\xi\in\bigoplus\limits_{\beta\in B} L_2(\mathbb{R}_\mu,\mu_\beta)
\]
means that
\[
\xi=\{\xi_\beta(\lambda)\}_{\beta\in B}, \quad \text{where} \quad
\xi_\beta: \mathbb{R}_\mu \to \mathbb{C}, \quad \text{and}\quad
\xi_\beta \in L_2(\mathbb{R}_\mu,\mu_\beta).
\]

Therefore,
\[
\int\limits_\mu^\infty |\xi_\beta(\lambda)|^2\,
d\mu_\beta(\lambda)\ <\ +\infty.
\]

Besides,
\begin{equation}\label{V}
(V\xi)_\beta(\lambda)\ =\ \lambda \,\xi_\beta(\lambda), \quad
\text{where} \quad V=UAU^{-1},
\end{equation}
that is, the projection of the operator $A$ onto the space
$L_2(\mathbb{R}_\mu,\mu_\beta)$ acts as a product by $\lambda$.

Let $w(t)$ be the solution of the homogeneous Cauchy problem
(\ref{prob2}). Our goal is to show that $w(t)=0$. Set
\[
\xi(t)\ =\ Uw(t).
\]

Here
\[
\xi(t)\ =\ \{\xi_\beta(t,\lambda)\}_{\beta\in B} , \quad
\text{where}\quad \xi_\beta:\mathbb{R}_0\times \mathbb{R}_\mu \to
\mathbb{C}.
\]

Then
\[
\partial_t^\rho \xi(t)\ =\ \partial_t^\rho Uw(t)\ =\ U\partial_t^\rho w(t)\ =\ UAw(t)\ =\
\]
\[
=\ UAU^{-1}Uw(t)\ =\ VUw(t)\ =\ V\xi(t).
\]

Further (see \cite{Dun}, Chapter XII, Sec. 3, Lemma 3 and
(\ref{V})),
\[
\partial_t^\rho \xi_\beta(t,\lambda)\ =\ (\partial_t^\rho\xi(t))_\beta(\lambda)\ =\
(V\xi(t))_\beta(\lambda)\ =\ \lambda\xi_\beta(t,\lambda).
\]

It is clear that
\[
\xi_\beta(0,\lambda) = 0.
\]

Thus, $\xi_\beta(t,\lambda)$ is the solution of the following
homogeneous Cauchy problem
\[
\partial_t^\rho \xi_\beta(t,\lambda)\ =\ \lambda\xi_\beta(t,\lambda),
\]
\[
\xi_\beta(0,\lambda) = 0.
\]

Consequently, $\xi_\beta(t,\lambda)\equiv 0$ for all $\beta$. Then
$\xi(t)=0$ and $w(t) = U^{-1}\xi(t)=0$.

\bibliographystyle{amsplain}

\end{document}